\documentclass[11pt]{article}

\usepackage{latexsym,amsmath,amssymb}

\begin{document}

\centerline{\large GENERALIZATION OF AN IDENTITY OF ANDREWS}

\vspace{1.5cm}

\centerline{\large Eduardo H.\ M.\ Brietzke}

\centerline{\small Instituto de Matem\'atica -- UFRGS}

\centerline{\small Caixa Postal 15080}

\centerline{\small 91509--900 Porto Alegre, RS, Brazil}

\centerline{\small email: brietzke@mat.ufrgs.br}

\vspace{0.7cm}

\begin{abstract}

We consider an identity relating Fibonacci numbers to Pascal's
triangle discovered by \mbox{G.\ E.\ Andrews}. Several authors
provided proofs of this identity, all of them rather involved or else
relying on sophisticated number theoretical arguments. We not only
give a simple and elementary proof, but also show the identity
generalizes to arrays other than Pascal's triangle. As an application
we obtain identities relating trinomial coefficients and Catalan's
triangle to Fibonacci numbers.

\end{abstract}

\vspace{2.3cm}

There is a well-known identity relating the sequence of
Fibonacci numbers to \mbox{Pascal's} triangle. Not so well-known are
identities
\begin{equation}
F_n=\sum_{k=-\infty}^\infty(-1)^k \,
{n-1 \choose \left[\frac{1}{2}\bigl(n-1-5k\bigr)\right]}
\hspace{0.15cm}
\end{equation}
and
\begin{equation}
F_n=\sum_{k=-\infty}^\infty(-1)^k \,
{n \choose \left[\frac{1}{2}\bigl(n-1-5k\bigr)\right]},
\end{equation}
obtained by G.\ E.\ Andrews in \cite{1}. Different proofs of (1) and
(2) were given by H.\ Gupta in \cite{4} and by M.\ D.\ Hirschhorn in
\cite{5} and \cite{6}. They are all
specifically
designed to deal with the case of Pascal's triangle. As indicated in
\cite{4} and \cite{5}, identities (1) and (2) are equivalent to
\begin{equation}
F_{2n+1}=\sum_{j=-\infty}^\infty
\biggl[{2n+1\choose n-5j}-{2n+1\choose n-5j-1}\biggr],
\end{equation}
\begin{equation}
F_{2n+2}=\sum_{j=-\infty}^\infty
\biggl[{2n+2\choose n-5j}-{2n+2\choose n-5j-1}\biggr]
\hspace{0.2cm}
\end{equation}
and
\begin{equation}
F_{2n+2}=\sum_{j=-\infty}^\infty
\biggl[{2n+1\choose n-5j}-{2n+1\choose n-5j-2}\biggr],
\end{equation}
\begin{equation}
F_{2n+1}=\sum_{j=-\infty}^\infty
\biggl[{2n\choose n-5j}-{2n\choose n-5j-2}\biggr],
\end{equation}
respectively. These identities have been reobtained by G.\
E.\ Andrews \cite{2} in the context of identities of the
Rogers--Ramanujan type (see also \cite{7}).

The purpose of this article is to provide a simple and entirely
elementary proof for identities (3) through (6) as well as several
other similar identities. Our method has the advantage of showing that
this kind of identity holds not only for Pascal's triangle, but also
for any other array constructed in a similar way. As an application
of our result, we  obtain versions of these identities for the
trinomial coefficients [the identities (9) and (10) below], as well
as for Catalan's triangle [(11) and (12)].

\vspace{0.3cm}

We start with a visualization of Pascal's triangle in which alternate
rows have been removed and only nonvanishing binomial numbers are
represented:

\hspace{2.cm}\begin{minipage}{9.4cm}

\begin{picture}(200,150)

\put(87,110){\mbox{\fbox{1}}}
\put(110,110){\mbox{1}}

\put(70,90){\mbox{1}}
\put(87,90){\mbox{\fbox{3}}}
\put(110,90){\mbox{3}}
\put(130,90){\mbox{1}}
\put(72.5,93.5){\circle{14}}

\put(50,70){\mbox{1}}
\put(70,70){\mbox{5}}
\put(85,70){\mbox{\fbox{10}}}
\put(107.5,70){\mbox{10}}
\put(130,70){\mbox{5}}
\put(150,70){\mbox{1}}
\put(72.5,73.5){\circle{14}}

\put(30,50){\mbox{1}}
\put(50,50){\mbox{7}}
\put(67.5,50){\mbox{21}}
\put(85,50){\mbox{\fbox{35}}}
\put(107.5,50){\mbox{35}}
\put(127.5,50){\mbox{21}}
\put(150,50){\mbox{7}}
\put(170,50){\mbox{1}}
\put(172.5,53.5){\circle{14}}
\put(72.5,53.5){\circle{15}}

\put(10,30){\mbox{1}}
\put(30,30){\mbox{9}}
\put(47.5,30){\mbox{36}}
\put(67.5,30){\mbox{84}}
\put(83.5,30){\mbox{\fbox{\!126\!}}}
\put(105.5,30){\mbox{126}}
\put(127.5,30){\mbox{84}}
\put(147.5,30){\mbox{36}}
\put(170,30){\mbox{9}}
\put(187,30){\mbox{\fbox{1}}}
\put(172.5,33.5){\circle{14}}
\put(72.5,33.5){\circle{15}}

\put(-10,15){\mbox{$\cdot$}}
\put(11,15){\mbox{$\cdot$}}
\put(31,15){\mbox{$\cdot$}}
\put(51,15){\mbox{$\cdot$}}
\put(71,15){\mbox{$\cdot$}}
\put(91,15){\mbox{$\cdot$}}
\put(111,15){\mbox{$\cdot$}}
\put(131,15){\mbox{$\cdot$}}
\put(151,15){\mbox{$\cdot$}}
\put(171,15){\mbox{$\cdot$}}
\put(191,15){\mbox{$\cdot$}}
\put(211,15){\mbox{$\cdot$}}

\end{picture}

\end{minipage}
\begin{minipage}{1cm}

(7)

\end{minipage}

\setcounter{equation}{7}

Identity (3) is visualised considering in each row the sum of the
elements represented inside rectangles and subtracting from it the
sum of elements represented inside circles. The remaining
identities can be likewise visualized. Note that the above array
consists of a first row in which the only nonvanishing elements are
two 1's followed by rows in which each element is obtained by adding
up the element above it to the right, the element above it to the
left, and 2 times the element directly above it. In general we have
the following theorem.

\vspace{0.5cm}

\noindent {\bf Theorem.} \emph{Let $s(0,k)$, with
$k\in\mathbb{Z}$, be an arbitrary sequence such that $s(0,k)\neq0$
for only finitely many values of $k$. Given
$\alpha,\beta\in\mathbb{R}$, define $s(n,k)$ recursively, for
$n\in\mathbb{N}$ and $k\in\mathbb{Z}$, setting
\[
s(n,k)=\alpha s(n-1,k-1)+\beta s(n-1,k)+\alpha s(n-1,k+1).
\]
Then for any fixed $k_0$ the sequence $\bigl(d_n\bigr)$
defined by
\[
d_n=\sum_{j=-\infty}^\infty
\Bigl[s(n,k_0-5j)-s(n,k_0-5j-1)\Bigr]
\]
satisfies the recurrence relation}
\[
d_n=\bigl(2\beta-\alpha\bigr)d_{n-1}+
\bigl(\alpha\beta+\alpha^2-\beta^2\bigr)d_{n-2}\,.
\]

\vspace{0.5cm}

\noindent {\bf Proof:} If suffices to prove the relation for $n=2$,
since the general case follows from it considering the row
$\bigl(s(n-2,k)\bigr)_k$ as being the first.

\vspace{0.2cm}

For each fixed $n$, the $n^{\rm th}$ row of the array depends linearly
on the first. In addition, the operator $\mathcal{L}$ defined by
\[
\mathcal{L}(a_n)=\sum_{j=-\infty}^\infty
\Bigl[a_{k_1-5j}-a_{k_1-5j-1}\Bigr]
\]
maps two-tailed sequences linearly to real numbers.

It therefore suffices to prove the proposition in the particular case
where the first row contains an element, say $s(0,0)$, equal to 1 and
all others equal to 0. Now it is a simple matter of inspecting the
table

\hspace{1.cm}
\begin{picture}(260,100)
\put(130,75){\mbox{1}}
\put(170,75){\mbox{0}}
\put(210,75){\mbox{0}}
\put(250,75){\mbox{0}}
\put(90,75){\mbox{0}}
\put(50,75){\mbox{0}}
\put(10,75){\mbox{0}}

\put(10,45){\mbox{0}}
\put(50,45){\mbox{0}}
\put(90,45){\mbox{$\alpha$}}
\put(130,45){\mbox{$\beta$}}
\put(170,45){\mbox{$\alpha$}}
\put(210,45){\mbox{0}}
\put(250,45){\mbox{0}}

\put(10,15){\mbox{0}}
\put(49,15){\mbox{$\alpha^2$}}
\put(82,15){\mbox{$2\alpha\beta$}}
\put(115,15){\mbox{$2\alpha^2\!+\!\beta^2$}}
\put(165,15){\mbox{$2\alpha\beta$}}
\put(209,15){\mbox{$\alpha^2$}}
\put(250,15){\mbox{0}}

\end{picture}

\noindent Since the null sequence satisfies any homogeneous linear
recurrence relation, it is enough to check that
$(d_0,d_1,d_2)=\pm(0,\alpha,2\alpha\beta-\alpha^2)$ and
$(d_0,d_1,d_2)=\pm(1,\beta-\alpha,2\alpha^2\!+\!\beta^2-2\alpha\beta)$
satisfy the relation $d_2=\bigl(2\beta-\alpha\bigr)d_1+
\bigl(\alpha\beta+\alpha^2-\beta^2\bigr)d_0$. This completes the
proof.

\vspace{0.5cm}

\noindent {\bf Corollary.} \emph{Under the same hypotheses as in the
theorem, for any given $k_0\in\mathbb{Z}$ and $k_1\in\mathbb{N}$,
the sequence $\bigl(d_n\bigr)$ defined by
\[
d_n=\sum_{j=-\infty}^\infty
\Bigl[s(n,k_0-5j)-s(n,k_0-5j-k_1)\Bigr]
\]
satisfies the recurrence relation
\[
d_n=\bigl(2\beta-\alpha\bigr)d_{n-1}+
\bigl(\alpha\beta+\alpha^2-\beta^2\bigr)d_{n-2}\,.
\]}


\noindent {\bf Proof:} The proof follows by linearity, writing out
\[
d_n=d_n^{(1)}+\cdots+d_n^{(k_1)}\,,
\]
with
\[
d_n^{(i)}=\sum_{j=-\infty}^\infty
\Bigl[s(n,k_0-i+1-5j)-s(n,k_0-i-5j)\Bigr],
\]
and applying the theorem.

\vspace{0.5cm}

\noindent {\bf Application 1.} Identities (3) through (6) hold. We
first apply  the theorem to the array (7), in which
$s(n,k)={2n+1\choose k+n}$, for $-n\le k\le n+1$, and 0 otherwise.
Since $\alpha=1$ and $\beta=2$, it follows that the right-hand side
of (3) defines a sequence satisfying the recurrence
relation $d_n=3d_{n-1}-d_{n-2}$. But it is a well-known fact that
this recurrence relation is also satisfied by both sequences
$d_n=F_{2n+1}$ and $d_n=F_{2n+2}$. Therefore, it suffices to verify
(3) for $n$ equal to 0 and 1. Identities (4) through (6) can be
obtained by the same argument, applied to (7) or to the array obtained
by deleting the even-numbered rows of Pascal's triangle.

\vspace{0.5cm}

\noindent {\bf Application 2.} The Fibonacci numbers appear in a
similar manner when we operate with the array of the trinomial
coefficients. The trinomial coefficients ${n\choose k}_{\!2}$, for
$|k|\le n$, are defined by (see \cite{3}, section 6.2)
\[
\bigl(1+x+x^2\bigr)^n=\sum_{k=-n}^n{n\choose k}_{\!\!2}x^{n+k}
\]
and satisfy
\[
{n\choose k}_{\!\!2}=\sum_j(-1)^j{n\choose j}{2n-2j\choose n-j-k}.
\]
Using the property
\begin{equation}
{n\choose k}_{\!\!2}={n-1\choose k-1}_{\!\!2}+
{n-1\choose k}_{\!\!2}+{n-1\choose k+1}_{\!\!2},
\end{equation}
we construct the array of the trinomial coefficients, as follows:

\hspace{-0.8cm}
\begin{picture}(370,150)

\put(9,2){\mbox{$\cdot$}}

\put(31,2){\mbox{$\cdot$}}
\put(61,2){\mbox{$\cdot$}}
\put(92,2){\mbox{$\cdot$}}
\put(122,2){\mbox{$\cdot$}}
\put(152,2){\mbox{$\cdot$}}
\put(182,2){\mbox{$\cdot$}}
\put(212,2){\mbox{$\cdot$}}
\put(242,2){\mbox{$\cdot$}}
\put(272,2){\mbox{$\cdot$}}
\put(302,2){\mbox{$\cdot$}}
\put(331,2){\mbox{$\cdot$}}

\put(353,2){\mbox{$\cdot$}}

\put(27,20){\mbox{\fbox{1}}}
\put(60,20){\mbox{5}}
\put(88,20){\mbox{15}}
\put(118,20){\mbox{30}}
\put(148,20){\mbox{45}}
\put(175,20){\mbox{\fbox{51}}}
\put(208,20){\mbox{45}}
\put(238,20){\mbox{30}}
\put(268,20){\mbox{15}}
\put(300,20){\mbox{5}}
\put(327,20){\mbox{\fbox{1}}}

\put(60,40){\mbox{1}}
\put(90,40){\mbox{4}}
\put(118,40){\mbox{10}}
\put(148,40){\mbox{16}}
\put(175,40){\mbox{\fbox{19}}}
\put(208,40){\mbox{16}}
\put(238,40){\mbox{10}}
\put(270,40){\mbox{4}}
\put(300,40){\mbox{1}}

\put(90,60){\mbox{1}}
\put(120,60){\mbox{3}}
\put(150,60){\mbox{6}}
\put(177,60){\mbox{\fbox{7}}}
\put(210,60){\mbox{6}}
\put(240,60){\mbox{3}}
\put(270,60){\mbox{1}}

\put(120,80){\mbox{1}}
\put(150,80){\mbox{2}}
\put(177,80){\mbox{\fbox{3}}}
\put(210,80){\mbox{2}}
\put(240,80){\mbox{1}}

\put(150,100){\mbox{1}}
\put(177,100){\mbox{\fbox{1}}}
\put(210,100){\mbox{1}}

\put(177,120){\mbox{\fbox{1}}}

\put(153,103){\circle{15}}
\put(153,83){\circle{15}}
\put(153,63){\circle{15}}
\put(154,43){\circle{15}}
\put(154,23){\circle{15}}

\put(303,23.5){\circle{15}}
\put(303,43.5){\circle{15}}

\end{picture}

\vspace{0.2cm}

We claim that the trinomial coefficients satisfy the identities
\begin{equation}
F_n=\sum_{j=-\infty}^\infty
\biggl[{n+1\choose 5j}_{\!\!2}-{n+1\choose 5j-1}_{\!\!2}\biggr]
\end{equation}
and
\begin{equation}
F_{n+1}=\sum_{j=-\infty}^\infty
\biggl[{n\choose 5j}_{\!\!2}-{n\choose 5j-2}_{\!\!2}\biggr].
\end{equation}

\noindent For example, identity (9) corresponds to adding in each row
the elements inside rectangles and subtracting the ones inside circles
in the above representation.

Since $\alpha=\beta=1$ by (8), our theorem and its corollary imply
that the right-hand sides of identities (9) and (10) define
sequences satisfying the recurrence relation $d_n=d_{n-1}+d_{n-2}$.
Hence, the proof of these identities consists in verifying them for
$n=0$ and 1.

\vspace{0.5cm}

\noindent {\bf Application 3.} On the bidimensional lattice
$\mathbb{Z}^2$, consider all paths that start at the origin, consist
of unit steps and are such that all steps go East or North. The length
of a path is the number of steps in the path. The distance between two
paths of length $n$ with end-points $(a_n,b_n)$ and $(a_n',b_n')$,
respectively, is
$|a_n-a_n'|$. Two paths are said to be non-intersecting if the origin
is the only point in common. Let
$B(n,k)$, for $1\le k\le n$, denote the number of pairs of
non-intersecting paths of length $n$ whose distance from one another
is $k$. The array defined by the $B(n,k)$ is called Catalan's triangle
(see \cite{8}) and its first column is formed by Catalan numbers
$C_n=B(n,1)=\frac{1}{n+1}{2n \choose n}$, which, among other things,
count rooted binary trees and the lattice paths from $(0,0)$ to
$(n,n)$ that stay above the main diagonal. Catalan's
triangle satisfies the same recurrence relation
\[
B(n,k)=B(n-1,k-1)+2B(n-1,k)+B(n-1,k+1)
\]
as alternate rows in Pascal's triangle, for $k\ge2$, while for $k=1$
we have
\[
B(n,1)=2B(n-1,1)+B(n-1,2).
\]

In order to fit Catalan's triangle into our framework, we embed it
into the larger array in which the first row is $s(1,\pm1)=\pm1$,
$s(1,k)=0$, if $k\neq\pm1$, and the subsequent rows are given by
\[
s(n,k)=s(n-1,k-1)+2s(n-1,k)+s(n-1,k+1), \quad
(n\ge2,\,k\in\mathbb{Z})
\]
as follows:

\hspace{-0.5cm}
\begin{picture}(360,145)

\put(5,2){\mbox{$\cdot$}}
\put(35,2){\mbox{$\cdot$}}
\put(65,2){\mbox{$\cdot$}}
\put(94,2){\mbox{$\cdot$}}
\put(123,2){\mbox{$\cdot$}}
\put(153,2){\mbox{$\cdot$}}
\put(180,2){\mbox{$\cdot$}}
\put(209,2){\mbox{$\cdot$}}
\put(236,2){\mbox{$\cdot$}}
\put(263,2){\mbox{$\cdot$}}
\put(290,2){\mbox{$\cdot$}}
\put(317,2){\mbox{$\cdot$}}
\put(343,2){\mbox{$\cdot$}}

\put(0,20){\mbox{$-1$}}
\put(26,20){\mbox{$-10$}}
\put(54,20){\mbox{$-44$}}
\put(82,20){\mbox{$-110$}}
\put(111,20){\mbox{$-165$}}
\put(140,20){\mbox{$-132$}}
\put(180,20){\mbox{0}}
\put(203,20){\mbox{132}}
\put(230,20){\mbox{165}}
\put(257,20){\mbox{110}}
\put(286,20){\mbox{44}}
\put(313,20){\mbox{10}}
\put(342,20){\mbox{1}}

\put(29,40){\mbox{$-1$}}
\put(56,40){\mbox{$-8$}}
\put(84,40){\mbox{$-27$}}
\put(113,40){\mbox{$-48$}}
\put(142,40){\mbox{$-42$}}
\put(180,40){\mbox{0}}
\put(205,40){\mbox{42}}
\put(232,40){\mbox{48}}
\put(259,40){\mbox{27}}
\put(288,40){\mbox{8}}
\put(315,40){\mbox{1}}

\put(57,60){\mbox{$-1$}}
\put(86,60){\mbox{$-6$}}
\put(113,60){\mbox{$-14$}}
\put(142,60){\mbox{$-14$}}
\put(180,60){\mbox{0}}
\put(205,60){\mbox{14}}
\put(232,60){\mbox{14}}
\put(261,60){\mbox{6}}
\put(288,60){\mbox{1}}

\put(86,80){\mbox{$-1$}}
\put(116,80){\mbox{$-4$}}
\put(145,80){\mbox{$-5$}}
\put(180,80){\mbox{0}}
\put(207,80){\mbox{5}}
\put(234,80){\mbox{4}}
\put(261,80){\mbox{1}}

\put(116,100){\mbox{$-1$}}
\put(145,100){\mbox{$-2$}}
\put(180,100){\mbox{0}}
\put(207,100){\mbox{2}}
\put(234,100){\mbox{1}}

\put(145,120){\mbox{$-1$}}
\put(180,120){\mbox{0}}
\put(207,120){\mbox{1}}

\end{picture}

\vspace{0.3cm}

\noindent The portion of the array lying to the the right of the zero
column is Catalan's triangle, i.e., $s(n,k)=B(n,k)$, if
$1\le k\le n$.

By the above theorem and corollary, several identities follow.
However, because $k\mapsto s(n,k)$ is odd, some of these identities
vanish trivially. We point out two nontrivial identities. Since
\[
\sum_{j=-\infty}^\infty\bigl[s(n,5j+1)-s(n,5j+4)\bigr]=
2\sum_{j=0}^\infty\bigl[B(n,5j+1)-B(n,5j+4)\bigr],
\]
using the same argument as before, we can show the identity
\begin{equation}
F_{2n-1}=\sum_{j=0}^\infty\bigl[B(n,5j+1)-B(n,5j+4)\bigr].
\end{equation}
Likewise one may show that
\begin{equation}
F_{2n-2}=\sum_{j=0}^\infty\bigl[B(n,5j+2)-B(n,5j+3)\bigr].
\end{equation}

According to \cite{8}, Proposition 2.1,
$B(n,k)=\frac{k}{n}{2n\choose n-k}$. Substituting in (13) and (14),
yields
\begin{equation}
F_{2n-1}=\sum_{j=0}^\infty\biggl[
\frac{5j+1}{n}{2n\choose n-5j-1}-\frac{5j+4}{n}{2n\choose n-5j-4}
\biggr]
\end{equation}
and
\begin{equation}
F_{2n-2}=\sum_{j=0}^\infty\biggl[
\frac{5j+2}{n}{2n\choose n-5j-2}-\frac{5j+3}{n}{2n\choose n-5j-3}
\biggr].
\end{equation}
Formulas (13) and (14) above appear to be new.

\vspace{0.5cm}

\noindent {\bf Acknowledgment.} The author would like to express
his gratitude to Prof.\ J.\ P.\ O.\ Santos for bringing identities
(3) and (5) to his attention. He is also grateful to Prof.\ J.\
Cigler for pointing out that a $q-$version of (1) and (2) already
appeared in the work \cite{27} by I.\ Schur from 1917.

\vspace{1.2cm}

\noindent AMS Classification Numbers: 11B39, 05A19

\end{document}